\documentclass[12pt,reqno]{amsart}
\usepackage{amsmath,amsrefs,amssymb}

\usepackage[colorlinks=true, pdfborder={0 0 0}]{hyperref}

\numberwithin{equation}{section}

\DeclareMathOperator{\divergence}{div}
\DeclareMathOperator{\spec}{spec}
\DeclareMathOperator{\range}{range}

\DeclareMathOperator{\smallo}{o}
\DeclareMathOperator{\Ricci}{Ric}
\newcommand{\R}{\mathbb{R}}
\renewcommand{\S}{\mathbb{S}}
\newcommand{\N}{\mathbb{N}}

\newcommand{\<}{\left<}
\renewcommand{\>}{\right>}
\renewcommand{\[}{\left[}
\renewcommand{\]}{\right]}
\renewcommand{\(}{\left(}
\renewcommand{\)}{\right)}

\newtheorem{theorem}{Theorem}[section]
\newtheorem{definition}[theorem]{Definition}

\begin{document}

\title[Non-synchronized solutions]{Non-synchronized solutions to nonlinear elliptic Schr\"odinger systems on a closed Riemannian manifold}

\author{Saikat Mazumdar}

\address{Saikat Mazumdar, Department of Mathematics, Indian Institute of Technology Bombay, Mumbai 400076, India}
\email{saikat@math.iitb.ac.in, saikat.mazumdar@iitb.ac.in}

\author{J\'er\^ome V\'etois}

\address{J\'er\^ome V\'etois, Department of Mathematics and Statistics, McGill University, 805 Sherbrooke Street West, Montreal, Quebec H3A 0B9, Canada}
\email{jerome.vetois@mcgill.ca}

\thanks{The second author was supported by the Discovery Grant RGPIN-2016-04195 from the Natural
Sciences and Engineering Research Council of Canada. This work was partially conducted when the first author held a postdoctoral position at McGill University under the co-supervision of Professors Pengfei Guan, Niky Kamran and the second author, that was supported by the NSERC Discovery Grants RGPIN-04443-2018, RGPIN-05490-2018 and RGPIN-04195-2016.}

\date{June 22, 2021}

\begin{abstract}
On a smooth, closed Riemannian manifold, we study the question of proportionality of components, also called synchronization, of vector-valued solutions to nonlinear elliptic Schr\"odinger systems with constant coefficients. In particular, we obtain bifurcation results showing the existence of branches of non-synchronized solutions emanating from the constant solutions.
\end{abstract}

\maketitle

\section{Introduction}\label{Intro}

On a smooth, closed Riemannian manifold $\(M,g\)$ of dimension $n$, we consider vector-valued solutions $\(u_1,u_2\)\in C^2\(M\)^2$ to elliptic systems of the form
\begin{equation}\label{IntroEq0}
\left\{\begin{aligned}&\Delta_{g}u_1=F_1\(u_1,u_2\)&&\text{in }M\\
&\Delta_{g}u_2=F_2\(u_1,u_2\),&&\text{in }M\end{aligned}\right.
\end{equation}
where $F_1$ and $F_2$ are $C^1$ functions and $\Delta_g:=-\divergence\nabla$ is the Laplace--Beltrami operator. In particular, we are interested in the stationary nonlinear Schr\"odinger system
\begin{equation}\label{IntroEq1}
\left\{\begin{aligned}&\Delta_g u_1+\lambda_1u_1=a_{11}u_1^{q-1}+a_{12}u_2^{q-2}u_1&&\text{in }M\\&\Delta_g u_2+\lambda_2u_2=a_{21}u_1^{q-2}u_2+a_{22}u_2^{q-1}&&\text{in }M\\&u_1,u_2>0&&\text{in }M,\end{aligned}\right.
\end{equation}
where $\lambda_1,\lambda_2,a_{11},a_{12},a_{21},a_{22}\in\R$ and $q\in\(2,\infty\)$. In the cubic case $q=4$, the system \eqref{IntroEq1} arises in particular in nonlinear optics (see for instance Akhmediev and Ankiewicz~\cite{AA} and Kanna and Lakshmanan~\cite{KL}) and the Hartree--Fock theory for Bose--Einstein condensates (see Esry, Greene, Burke and Bohn~\cite{EGBB} and Timmermans~\cite{T}). Such systems have received considerable attention from mathematicians in recent years. Among many others, let us refer for instance to the work of Bartsch, Dancer and Wang~\cite{BDW}, Clapp and Pistoia~\cites{CP1,CP2}, Gladiali, Grossi and Troestler~\cites{GGT2,GGT1}, Guo, Li and Wei~\cite{GLW}, Guo and Liu~\cite{GL}, Lin and Wei~\cite{LW1}, Liu and Wang~\cite{LW2}, Peng, Peng and Wang~\cite{PPW}, Sirakov~\cite{S}, Soave and Zilio~\cite{SZ}, Terracini and Verzini~\cite{TV} and Wei and Wu~\cite{WW} in the case where $M=\R^n$ (note that when $\(M,g\)$ is the standard round sphere, $n\ge3$, $\lambda_1=\lambda_2=n\(n-2\)/4$ and $q=2n/\(n-2\)$, we can use stereographic projection to write \eqref{IntroEq1} as a system in $\R^n$) and Chen and Zou~\cite{CZ}, Clapp, Pistoia and Tavares~\cite{CPT}, Druet and Hebey~\cite{DH} and Druet, Hebey and V\'etois~\cite{DHV} in the case of a more general manifold.

\smallskip
In this paper, we are interested in the question of proportionality of components, also called synchronization, of solutions to the system \eqref{IntroEq1}. A solution $\(u_1,u_2\)$ of \eqref{IntroEq1} is said to be {\it synchronized} if there exists a constant $\Lambda>0$ such that $u_2\equiv\Lambda u_1$ in $M$. This question has been studied for instance by Montaru, Sirakov and Souplet~\cite{MSS} and Quittner and Souplet~\cite{QS} in the case of systems in domains of $\R^n$. It also naturally arises in the case of a closed manifold.

\smallskip
It is easy to see that every synchronized solution $\(u_1,u_2\)=\(u_1,\Lambda u_1\)$ of the system \eqref{IntroEq1} is such that $u_1$ is constant in the case where $\lambda_1\ne\lambda_2$ and $u_1$ is a solution of the equation
\begin{equation}\label{IntroEq2}
\left\{\begin{aligned}&\Delta_g u_1+\lambda_1 u_1=\mu_1 u_1^{q-2}&&\text{ in }M\\&u_1>0&&\text{ in }M,\end{aligned}\right.
\end{equation}
where $\mu_1:=a_{11}+a_{12}\Lambda^{q-2}=a_{21}+a_{22}\Lambda^{q-2}$, in the case where $\lambda_1=\lambda_2$. We know from a result of Bidaut-V\'eron and V\'eron~\cite{BV} that the equation \eqref{IntroEq2} does not have any non-constant solutions when 
\begin{equation}\label{IntroEq3}
2<q\le2^*\quad\text{and}\quad\left\{\begin{aligned}&\frac{n-1}{n}\(q-2\)\lambda_1g\le\Ricci_g&&\text{if }q<2^*\\&\frac{n-1}{n}\(q-2\)\lambda_1g<\Ricci_g&&\text{if }q=2^*,\end{aligned}\right.
\end{equation}
where $2^*:=\infty$ if $n\le2$, $2^*:=2n/\(n-2\)$ if $n\ge3$, $\Ricci_g$ is the Ricci curvature of the manifold and the latter inequalities are in the sense of bilinear forms. On the other hand, existence results of non-constant solutions to equation \eqref{IntroEq2} abound in the case where \eqref{IntroEq3} is not satisfied (see for instance Chen, Wei and Yan~\cite{CWY}, Hebey and Vaugon~\cite{HV}, Hebey and Wei~\cite{HW}, Micheletti, Pistoia and V\'etois~\cite{MPV} and V\'etois and Wang~\cite{VW}.)

\smallskip
For simplicity, in this introduction, we state our results in the case of the sphere $\(\S^n,g_0\)$, where $g_0$ is the standard round metric. Furthermore, we assume that $\lambda_1=\lambda_2$ and $a_{12}=a_{21}$, namely we consider the system
\begin{equation}\label{IntroEq4}
\left\{\begin{aligned}&\Delta_{g_0} u_1+\lambda u_1=au_1^{q-1}+bu_2^{q-2}u_1&&\text{in }\S^n\\&\Delta_{g_0} u_2+\lambda u_2=bu_1^{q-2}u_2+cu_2^{q-1}&&\text{in }\S^n\\&u_1,u_2>0&&\text{in }\S^n,\end{aligned}\right.
\end{equation}
where $\lambda,a,b,c\in\R$ and $q\in\(2,\infty\)$. In the Euclidean space, this case has been studied for instance by Clapp and Pistoia~\cite{CP1} (via stereographic projection, the system studied in~\cite{CP1} matches with \eqref{IntroEq4} when $n=4$, $\lambda=2$ and the parameters $\alpha$ and $\beta$ in~\cite{CP1} are equal to 2). We refer to Sections~\ref{Bif} and~\ref{Proportional} for results applying to more general systems and more general manifolds. For the system~\eqref{IntroEq4}, we obtain the following:

\begin{theorem}\label{Th1}
Let $\lambda,a,b,c\in\R$ and $q\in\(2,\infty\)$. 
\begin{enumerate}
\item[(i)] If either $c\le b\le a$ or $a\le b\le c$ and at least one of the two inequalities is strict, then the system \eqref{IntroEq4} has no solutions.
\item[(ii)] If either [$a<b$ and $c<b$] or $a=b=c$, then every solution of \eqref{IntroEq4} is synchronized.
\item[(iii)] There exist real numbers $\lambda$, $a$, $b$ and $c$ such that $\lambda>0$, $a=c>b>0$ and \eqref{IntroEq4} has non-synchronized solutions. More precisely, we have the following result: for every $\lambda,a,b\in C^1\(\[-\delta,\delta\]\)$, $\delta>0$, if the following conditions hold:
\begin{enumerate}
\item[(A1)]$\lambda\(0\)\(a\(0\)+b\(0\)\)>0$,
\item[(A2)]$\lambda\(0\)\not\in\big\{\frac{2j\(2j+n-1\)}{q-2}:\ j\in\N\big\}$, where $\N:=\left\{1,2,\dotsc\right\}$,
\item[(A3)]$\beta\(0\):=\lambda\(0\)\frac{a\(0\)-b\(0\)}{a\(0\)+b\(0\)}\in\big\{\frac{j\(j+n-1\)}{q-2}:\ j\in\N\big\}$ and $\beta'\(0\)\ne0$,
\end{enumerate}
then there exists a $C^1$ branch (see Definition~\ref{Def}) of non-syn\-chronized solutions to \eqref{IntroEq4} with $\lambda=\lambda\(\alpha\)$, $c=a=a\(\alpha\)$ and $b=b\(\alpha\)$ emanating from the constant solution at $\alpha=0$.
\end{enumerate}
\end{theorem}

Theorem~\ref{Th1}~(iii) extends a previous result obtained by Gladiali, Grossi and Troestler~\cite{GGT1} for systems with Sobolev critical growth in $\R^n$, which, via stereographic projection, corresponds to the case where $n\ge3$, $\lambda=n\(n-2\)/4$ and $q=2n/\(n-2\)$. Like in~\cite{GGT1}, our approach is based on the bifurcation theory at eigenvalues of odd multiplicity. Unlike in~\cite{GGT1}, by taking advantage of our closed manifold setting, we perform our constructions in $C^{1,\theta}\(M\)$, $\theta\in\(0,1\)$, instead of Sobolev spaces, which allows us to treat the case of systems with supercritical growth. 

\smallskip
Theorem~\ref{Th1}~(iii) is proven in Section~\ref{Bif} (as a particular case of Theorem~\ref{Th3}) and Theorem~\ref{Th1}~(i) and~(ii) are proven in Section~\ref{Proportional} (as a particular cases of Theorem~\ref{Th4}~(i) and~(ii)).

\section{Bifurcation results}\label{Bif} 

This section is devoted to bifurcation results showing the existence of branches of non-synchronized solutions for systems like \eqref{IntroEq1}. 

\begin{definition}\label{Def}
Let $\(M,g\)$ be a smooth, closed Riemannian manifold, $\Omega$ be an open set in $\R^2$, $I:=\[-\delta,\delta\]$, $\delta>0$, and $F_1,F_2\in C^1\(I\times\Omega\)$. Consider the system
\begin{equation}\label{DefEq}
\left\{\begin{aligned}&\Delta_{g}u_1=F_1\(\alpha,u_1,u_2\)&&\text{in }M\\
&\Delta_{g}u_2=F_2\(\alpha,u_1,u_2\)&&\text{in }M,\end{aligned}\right.
\end{equation}
where $\alpha\in I$. Assume that for every $\alpha\in I$, there exists a solution $\(\overline{u}_{1}\(\alpha\), \overline{u}_{2}\(\alpha\)\)\in C^2\(M\)^2$ of \eqref{DefEq} such that $\(\overline{u}_{1}\(\alpha\), \overline{u}_{2}\(\alpha\)\)\to\(\overline{u}_{1}\(0\), \overline{u}_{2}\(0\)\)$ in $C^2\(M\)^2$ as $\alpha\to0$. Let $\mathcal{S}$ be the set of all solutions $\(\alpha,u_1, u_2\)\in I\times C^2\(M\)^2$ to \eqref{DefEq} such that $\(u_{1},u_{2}\)\ne\(\overline{u}_{1}\(\alpha\), \overline{u}_{2}\(\alpha\)\)$. We say that the solution $\(0,\overline{u}_{1}\(0\), \overline{u}_{2}\(0\)\)$ is a \emph{bifurcation point} of \eqref{DefEq} if $\(0,\overline{u}_{1}\(0\), \overline{u}_{2}\(0\)\)\in\overline{\mathcal{S}}$, where $\overline{\mathcal{S}}$ stands for the closure of $\mathcal{S}$ in $I\times C^2\(M\)^2$. Furthermore, we say that a subset $\mathcal{B}\subseteq\mathcal{S}$ is a $C^1$ \emph{branch} of solutions to \eqref{DefEq} emanating from $\(0,\overline{u}_{1}\(0\), \overline{u}_{2}\(0\)\)$ if $\mathcal{B}\ne\emptyset$ and $\mathcal{B}\cup\left\{\(0,\overline{u}_{1}\(0\),\overline{u}_{2}\(0\)\)\right\}$ is a $C^1$ curve in $I\times C^2\(M\)^2$.
\end{definition}

In the case of the sphere, we obtain Theorem~\ref{Th1}~(iii). In the case of a more general manifold, we obtain the following:

\begin{theorem}\label{Th2}
Let $\(M,g\)$ be a smooth, closed Riemannian manifold, $I:=\[-\delta,\delta\]$, $\delta>0$, $\lambda_1,\lambda_2,a_{11},a_{12},a_{21},a_{22}\in C^1\(I\)$ and $q\in\(2,\infty\)$. Consider the system
\begin{equation}\label{Th2Eq}
\left\{\begin{aligned}&\Delta_g u_1+\lambda_1\(\alpha\)u_1=a_{11}\(\alpha\)u_1^{q-1}+a_{12}\(\alpha\)u_2^{q-2}u_1&&\text{in }M\\&\Delta_g u_2+\lambda_2\(\alpha\)u_2=a_{21}\(\alpha\)u_1^{q-2}u_2+a_{22}\(\alpha\)u_2^{q-1}&&\text{in }M\\&u_1,u_2>0&&\text{in }M,\end{aligned}\right.
\end{equation}
where $\alpha\in I$. Assume that the following conditions hold:
\begin{enumerate}
\item[(B1)]
For every $\alpha \in I$, the system \eqref{IntroEq1} has a unique constant solution $\(\overline{u}_{1}\(\alpha\), \overline{u}_{2}\(\alpha\)\)$.
\item[(B2)]
For every $\alpha \in I$, the matrix 
\begin{align*}
\mathcal{A}\(\alpha\):= \(\begin{array}{cc} a_{11}\(\alpha\)\overline{u}_{1}\(\alpha\)^{q-2}& a_{12}\(\alpha\)\overline{u}_{2}\(\alpha\)^{q-3}\overline{u}_{1}\(\alpha\)\\a_{21}\(\alpha\) \overline{u}_{1}\(\alpha\)^{q-3}\overline{u}_{2}\(\alpha\)  & a_{22}\(\alpha\)\overline{u}_{2}\(\alpha\)^{q-2}\end{array}\)
\end{align*}
has two distinct, non-zero, real eigenvalues $\beta_{1}\(\alpha\)$ and $\beta_{2}\(\alpha\)$.
\item[(B3)] $\mathcal{H}^{*}_1\times\mathcal{H}^{*}_2$ has odd dimension, where
$$\mathcal{H}^{*}_{i}:=\left\{ \varphi \in C^{2}\(M\):\Delta_{g}\varphi=\(q-2\)\beta_{i}\(0\)\varphi\hbox{ in }M\right\}\quad\forall i\in\left\{1,2\right\}.$$
\item[(B4)] For every $i\in\left\{1,2\right\}$, if $\mathcal{H}^*_i\ne\left\{0\right\}$, then $\beta'_{i}\(0\)\ne 0$ and either $\lambda_1\(0\)\ne\lambda_2\(0\)$ or $\beta_{i}\(0\)\ne\lambda_1\(0\)=\lambda_2\(0\)$.
\end{enumerate}
Then the solution $\(0,\overline{u}_{1}\(0\), \overline{u}_{2}\(0\)\)$ is a bifurcating point of the system \eqref{Th2Eq}. Furthermore, there exists a neighborhood $\mathcal{N}$ of $\(0,\overline{u}_{1}\(0\), \overline{u}_{2}\(0\)\)$ in $I\times C^2\(M\)^2$ such that for every solution $\(\alpha,u_1, u_2\)\in\mathcal{N}$ of \eqref{Th2Eq}, if $\(u_{1},u_{2}\)\ne\(\overline{u}_{1}\(\alpha\), \overline{u}_{2}\(\alpha\)\)$, then $\(u_{1},u_{2}\)$ is non-synchronized. If moreover $\mathcal{H}^{*}_1\times\mathcal{H}^{*}_2$ has dimension one, then there exists a $C^1$ branch of non-synchronized solutions to \eqref{Th2Eq} emanating from $\(0,\overline{u}_{1}\(0\), \overline{u}_{2}\(0\)\)$.
\end{theorem}

Both Theorem~\ref{Th1}~(iii) and Theorem~\ref{Th2} follow from the following general bifurcation result for systems of the form \eqref{DefEq}:

\begin{theorem}\label{Th3}
Let $\(M,g\)$ be a smooth, closed Riemannian manifold, $\Omega$ be an open set in $\R^2$ such that $\(0,0\)\in\Omega$, $I:=\[-\delta,\delta\]$, $\delta>0$, and $F_1,F_2\in C^1\(I\times\Omega\)$ such that $\partial_\alpha\partial_{u_{i}}F_{j}$ exists and is continuous in $I\times\Omega$ for all $i,j\in\left\{1,2\right\}$, where we denote by $\(\alpha,u_1,u_2\)$ a point in $I\times\Omega$. Assume that the following conditions hold:
\begin{enumerate}
\item[(C1)] $F_1\(\alpha,0,0\)=F_2\(\alpha,0,0\)=0$  for all $\alpha\in I$. 
\item[(C2)]  $\partial _{u_2}F_1\(\alpha,0,0\)=\partial_{u_1} F_2\(\alpha,0,0\)=0$ for all $\alpha \in I$.
\item[(C3)] There exist two closed subspaces $\mathcal{H}_{1}$ and $\mathcal{H}_{2}$ of $C^{1,\theta}\(M\)$, $\theta\in\(0,1\)$, and two open subsets $\mathcal{U}_1\subseteq\mathcal{H}_1$ and $\mathcal{U}_2\subseteq\mathcal{H}_2$ which contain 0 and satisfy the following conditions:
\begin{itemize}
\item[--]$\(u_{1}\(x\),u_{2}\(x\)\)\in\Omega$ for all $x\in M$ and $\(u_1,u_2\)\in{\mathcal{U}_1}\times{\mathcal{U}_2}$.
\item[--]$ \(\Delta_{g}+1\)^{-1}u_{i}\in\mathcal{H}_{i}$ for all $u_{i}\in {{\mathcal{U}}_{i}}$ and $i\in\left\{1,2\right\}$.
\item[--]$\(\Delta_{g}+1\)^{-1}F_{i}\(\alpha,u_{1},u_{2}\) \in\mathcal{H}_{i}$ for all $\(\alpha,u_{1},u_{2}\)\in I\times{{\mathcal{U}}_{1}}\times{{\mathcal{U}}_{2}}$ and $i\in\left\{1,2\right\}$. 
\end{itemize} 
\item[(C4)] $\mathcal{H}^{*}_{1}\times \mathcal{H}^{*}_{2}$ has odd dimension, where 
$$\mathcal{H}^{*}_{i}:=\left\{\varphi\in \mathcal{H}_{i}:\ \Delta_{g}\varphi=\partial_{u_{i}}F_{i}\(0,0,0\)\varphi\text{ in }M\right\}\quad\forall i\in\left\{1,2\right\}.$$
\item[(C5)]For every $i\in\left\{1,2\right\}$, if $\mathcal{H}^{*}_{i}\ne\left\{0\right\}$, then $\partial_\alpha\partial_{u_{i}} F_{i}\(0,0,0\)\ne0$.
\end{enumerate}
Then there exists a sequence of solutions $\(\(\alpha_m,u_{1,m},u_{2,m}\)\)_{m\in\N}$ to the system \eqref{DefEq} such that $\(\alpha_m,u_{1,m},u_{2,m}\)\in I\times\(\(\mathcal{U}_{1}\times \mathcal{U}_{2}\)\backslash\left\{\(0,0\)\right\}\)$ and $\(\alpha_m,u_{1,m},u_{2,m}\)\to\(0,0,0\)$ in $I\times C^2\(M\)^2$ as $m\to\infty$. Furthermore, every such sequence $\(\(\alpha_m,u_{1,m},u_{2,m}\)\)_{m\in\N}$ is such that up to a subsequence,
\begin{equation}\label{bif:expan}
u_{i,m}=\varepsilon_{m}\varphi_{i}+\smallo\(\varepsilon_{m}\)\quad\hbox{in } \mathcal{H}_{i}\quad \forall i\in\left\{1,2\right\}
\end{equation}
as $m\to\infty$ for some $\(\varphi_{1},\varphi_{2}\)\in\(\mathcal{H}^{*}_{1}\times\mathcal{H}^{*}_{2}\)\backslash\left\{\(0,0\)\right\}$ and $\varepsilon_{m}>0$ such that $\varepsilon_m\to0$. If moreover $\mathcal{H}^{*}_1\times\mathcal{H}^{*}_2$ has dimension one and $\partial_{u_{i}}\partial_{u_{j}}F_{k}$ exists and is continuous in $I\times\Omega$ for all $i,j,k\in\left\{1,2\right\}$, then there exists a neighborhood $\mathcal{N}$ of $\(0,0,0\)$ in $I\times \mathcal{U}_{1}\times\mathcal{U}_{2}$ such that the set of solutions $\(\alpha,u_1,u_2\)\in\mathcal{N}\backslash\(I\times\left\{\(0,0\)\right\}\)$ to \eqref{DefEq} is a $C^1$ branch emanating from $\(0,0,0\)$ whose tangent line at $\(0,0,0\)$ is directed by some vector in $\R\times\(\(\mathcal{H}^{*}_{1}\times\mathcal{H}^{*}_{2}\)\backslash\left\{\(0,0\)\right\}\)$.
\end{theorem}

\proof[Proof of Theorem~\ref{Th3}]
Let 
$$\mathcal{H}:=\mathcal{H}_{1}\times\mathcal{H}_{2},\quad\mathcal{H}^*:=\mathcal{H}_{1}^*\times\mathcal{H}_{2}^*\quad\hbox{and}\quad\mathcal{U}:=\mathcal{U}_{1}\times\mathcal{U}_{2}.$$
By replacing $\mathcal{U}_{1}$ and $\mathcal{U}_{2}$ by smaller sets if necessary, we may assume that there exists a compact set $K\subset\Omega$ such that $\(u_{1}\(x\),u_{2}\(x\)\)\in K$ for all $x\in M$ and $\(u_1,u_2\)\in\mathcal{U}$. The solutions $\(u_1,u_2\)\in\mathcal{U}$ to \eqref{DefEq} are given by  the zeros of the function $\mathcal{T}:I\times{{\mathcal{U}}}\to\mathcal{H}$ defined by
$$\mathcal{T}\(\alpha,u_{1},u_{2}\):=\mathcal{I}\(u_1,u_2\)-\mathcal{K}_{\alpha}\(u_{1},u_{2}\),\quad\text{where}\quad\mathcal{I}\(u_1,u_2\):=\(\begin{array}{c}u_{1}\\u_{2}\end{array}\)$$
and 
$$\mathcal{K}_{\alpha}\(u_{1},u_{2}\):=\(\begin{array}{c}\(\Delta_{g}+1\)^{-1}\[F_{1}\(\alpha,u_{1},u_{2}\)+u_{1}\]\\\(\Delta_{g}+1\)^{-1}\[F_{2}\(\alpha,u_{1},u_{2}\)+u_{2}\]\end{array}\)$$
for all $\(\alpha,u_1,u_2\)\in I\times{{\mathcal{U}}}$. By the assumption (C3), we have that the functions $\mathcal{K}_{\alpha}$ and $\mathcal{T}$ are well-defined. In what follows, we write a point in $I\times{{\mathcal{U}}}$ as $\(\alpha,U\)$.

\smallskip\noindent {\it Step 1.}  
We begin with proving that $\mathcal{K}_{\alpha}$ is compact, $\mathcal{T}\in C^{1}\(I\times \mathcal{U}\)$ and $D_{U}\,\partial_\alpha\mathcal{T}$ exists and is continuous in $I\times \mathcal{U}$. 

\smallskip
Suppose that $\(\(u_{1,m},u_{2,m}\)\)_{m\in\N}$ is a bounded sequence in ${{\mathcal{U}}}$. Then the sequences $\(u_{1,m}\)_m$ and $\(u_{2,m}\)_m$ are bounded in $C^{1,\theta}\(M,K\)$ and up to a subsequence, $\(u_{1,m},u_{2,m}\)\to\(u_{1,0},u_{2,0}\)$ in $C^1\(M,K\)^2$ as $m\to\infty$. Let $\alpha\in I$ and 
$$\(\begin{array}{c}\tilde{u}_{1,m}\\\tilde{u}_{2,m}\end{array}\):=\mathcal{K}_{\alpha}\(u_{1,m},u_{2,m}\)$$
so that
$$\(\begin{array}{c}\Delta_{g}\tilde{u}_{1,m}+\tilde{u}_{1,m}\\\Delta_{g}\tilde{u}_{2,m}+ \tilde{u}_{2,m}\end{array}\)=\(\begin{array}{c}F_{1}\(\alpha,u_{1,m},u_{2,m}\)+u_{1,m}\\  F_{2}\(\alpha,u_{1,m},u_{2,m}\)+u_{2,m}\end{array}\).$$
Since $F_1,F_2\in C^1\(I\times\Omega\)$, $u_{1,0},u_{2,0}\in C^1\(M,K\)$ and $K\subset\Omega$, by standard elliptic estimates, we obtain that  there exists $\(\tilde{u}_{1,0},\tilde{u}_{2,0}\)\in C^{2}\(M\)^2$ satisfying
$$\(\begin{array}{c}\Delta_{g}\tilde{u}_{1,0}+\tilde{u}_{1,0}\\\Delta_{g}\tilde{u}_{2,0}+\tilde{u}_{2,0}\end{array}\)=\(\begin{array}{c}F_{1}\(\alpha,u_{1,0},u_{2,0}\)+u_{1,0}\\F_{2}\(\alpha,u_{1,0},u_{2,0}\)+u_{2,0}\end{array}\).$$
For $i=1,2$, we then obtain
\begin{multline}\label{arg_con}
\(\Delta_{g}+1\)\(\tilde{u}_{i,m}-\tilde{u}_{i,0}\)=F_{i}\(\alpha,u_{1,m},u_{2,m}\)-F_{i}\(\alpha,u_{1,0},u_{2,0}\)\\
+u_{i,m}-u_{i,0}=\smallo\(1\)
\end{multline}
uniformly in $M$. Therefore, $\(\tilde{u}_{1,m},\tilde{u}_{2,m}\)\to\(\tilde{u}_{1,0},\tilde{u}_{2,0}\)$ in $\mathcal{H}$ and so $\mathcal{K}_{\alpha}$ is compact.

\smallskip
Now, if $\(\alpha_{m},u_{1,m},u_{2,m}\)\to\(\alpha_{0},u_{1,0},u_{2,0}\)$ in $I\times{\mathcal{U}}$, then, arguing as above, we obtain that $\mathcal{K}_{\alpha_m}\(u_{1,m},u_{2,m}\)\to\mathcal{K}_{\alpha_0}\(u_{1,0},u_{2,0}\)$ in $\mathcal{H}$.
This shows that  $\mathcal{T}$ is continuous in $I\times{\mathcal{U}}$.

\smallskip
For every $\(\alpha,u_{1},u_{2}\)\in I\times{\mathcal{U}}$, we have
\begin{align*}
&\partial_\alpha\mathcal{T}\(\alpha,u_{1},u_{2}\)=-\(\begin{array}{c}\(\Delta_{g}+1\)^{-1}\[\partial_\alpha F_{1}\(\alpha,u_{1},u_{2}\)\]\\\(\Delta_{g}+1\)^{-1}\[\partial_\alpha F_{2}\(\alpha,u_{1},u_{2}\)\]\end{array}\),\allowdisplaybreaks\\
&D_{U}\mathcal{T}\(\alpha,u_{1},u_{2}\)\[\(v_{1},v_{2}\)\]\\
&=\(\begin{array}{c}v_1-\(\Delta_{g}+1\)^{-1}\[\partial_{u_{1}}F_{1}\(\alpha,u_{1},u_{2}\)v_{1}+\partial_{u_{2}}F_{1}\(\alpha,u_{1},u_{2}\)v_{2}+v_{1}\]\\v_2-\(\Delta_{g}+1\)^{-1}\[\partial_{u_{1}}F_{2}\(\alpha,u_{1},u_{2}\)v_{1}+\partial_{u_{2}}F_{2}\(\alpha,u_{1},u_{2}\)v_{2}+v_{2}\]\end{array}\)
\end{align*}
and
\begin{align*}
&D_U\,\partial_\alpha\mathcal{T}\(\alpha,u_{1},u_{2}\)\[\(v_{1},v_{2}\)\]\\
&\quad=-\(\begin{array}{c}\(\Delta_{g}+1\)^{-1}\[\partial_{u_{1}}\partial_{\alpha}F_{1}\(\alpha,u_{1},u_{2}\)v_{1}+\partial_{u_{2}}\partial_{\alpha}F_{1}\(\alpha,u_{1},u_{2}\)v_{2}\]\\\(\Delta_{g}+1\)^{-1}\[\partial_{u_{1}}\partial_{\alpha}F_{2}\(\alpha,u_{1},u_{2}\)v_{1}+\partial_{u_{2}}\partial_{\alpha}F_{2}\(\alpha,u_{1},u_{2}\)v_{2}\]\end{array}\).
\end{align*}
If $\(\alpha_{m},u_{1,m},u_{2,m}\)\to\(\alpha_{0},u_{1,0},u_{2,0}\)$ in $I\times{\mathcal{U}}$, then for $i,j=1,2$, by using the regularity assumptions on $F$, we obtain that
\begin{align*}
\partial_\alpha F_j\(\alpha_{m},u_{1,m},u_{2,m}\)&\longrightarrow\partial_\alpha F_j\(\alpha_0,u_{1,0},u_{2,0}\),\allowdisplaybreaks\\
\partial_{u_i} F_j\(\alpha_{m},u_{1,m},u_{2,m}\)&\longrightarrow \partial_{u_i}F_j\(\alpha_0,u_{1,0},u_{2,0}\)
\end{align*}
and 
$$\partial_{u_i}\partial_\alpha F_j\(\alpha_{m},u_{1,m},u_{2,m}\)\longrightarrow \partial_{u_i}\partial_\alpha F_j\(\alpha_0,u_{1,0},u_{2,0}\)$$
uniformly in $M$. Then, arguing as in \eqref{arg_con}, we obtain that $\mathcal{T}\in C^{1}\(I\times\mathcal{U}\)$ and $D_U\,\partial_\alpha\mathcal{T}$ is continuous in $I\times\mathcal{U}$.

\smallskip\noindent {\it Step 2.} 
We now establish the main bifurcation results.

\smallskip
We have $\mathcal{T}\(\alpha,0,0\)=\(0,0\)$ for all $\alpha\in I$ by the assumption (C1). Furthermore, by the assumption (C2), we have
\begin{align*}
&D_{U}\mathcal{T}\(0,0,0\)\[\(v_{1}, v_{2}\)\]=\(\begin{array}{c}v_1-\(\Delta_{g}+1\)^{-1}\[\partial_{u_{1}}F_{1}\(0,0,0\)v_{1}+v_{1}\]\\v_2-\(\Delta_{g}+1\)^{-1}\[\partial_{u_{2}}F_{2}\(0,0,0\)v_{2}+v_{2}\]\end{array}\).
\end{align*}
Then $\(v_{1},v_{2}\)\in\ker D_{U}\mathcal{T}\(0,0,0\)$ if and only if $\(v_{1},v_{2}\)\in\mathcal{H}^*$. By the assumption (C4), it follows that $D_{U}\mathcal{T}\(0,0,0\) $ has a nontrivial kernel consisting of eigenfunctions of $\Delta_{g}$.

\smallskip
Now, $\ker D_{U}\mathcal{T}\(0,0,0\)$ has odd dimension by the assumption (C4). This along with the results of Step~1 allows to apply Theorem~A of Westreich~\cite{W} (see also Theorems~II.3.3 and~II.4.4 and the statements (II.4.29) and (II.4.31) in Kielh\"{o}fer's book~\cite{K}), which gives that the solution $\(0,0,0\)$ is a bifurcation point of the system \eqref{DefEq} in $I\times\(\mathcal{H}\backslash\left\{\(0,0\)\right\}\)$ provided the following condition holds:
\begin{multline}\label{odd_bif}
\Big[D_U\partial_\alpha\mathcal{T}\(0,0,0\)\[\(v_{1},v_{2}\)\]\in\range\(D_{U}\mathcal{T}\(0,0,0\)\)\\
\text{and}\quad\(v_{1},v_{2}\)\in\ker\(D_{U}\mathcal{T}\(0,0,0\)\)\Big]\iff\(v_{1},v_{2}\)=\(0,0\).
\end{multline}
Remark that by standard elliptic estimates, the $C^2$ topology in Definition~\ref{Def} can be replaced without loss of generality by the $C^{1,\theta}$ topology. If moreover $\mathcal{H}^*$ has dimension one and $\partial_{u_{i}}\partial_{u_{j}}F_{k}$ exists and is continuous in $I\times\Omega$ for all $i,j,k\in\left\{1,2\right\}$, then the last part of Theorem~\ref{Th3} follows from Theorem~1.7 of Crandall--Rabinowitz~\cite{CR} (see also Kielh\"{o}fer~\cite{K}*{Theorem~I.5.1}).

\smallskip
We now show that the condition \eqref{odd_bif} holds for the function $\mathcal{T}$ under our assumptions on $F$. Let $\(v_{1},v_{2}\)\in\ker D_{U}\mathcal{T}\(0,0,0\)$ and $\(w_{1},w_{2}\):=D_U\partial_\alpha\mathcal{T}\(0,0,0\)\[\(v_{1},v_{2}\)\]$. Then 
$$\(\begin{array}{c}\Delta_{g}v_{1}\\\Delta_{g}v_{2}\end{array}\)=\(\begin{array}{c}\partial_{u_{1}}F_{1}\(0,0,0\)v_{1}\\\partial_{u_{2}}F_{2}\(0,0,0\)v_{2}\end{array}\).$$
Furthermore, by the assumption (C2), we obtain
$$\(\begin{array}{c}\Delta_{g}w_{1}+w_{1}\\\Delta_{g}w_{2}+w_{2}\end{array}\)=-\(\begin{array}{c}\partial_{u_{1}}\partial_{\alpha}F_{1}\(0,0,0\)v_{1}\\\partial_{u_{2}}\partial_{\alpha}F_{2}\(0,0,0\)v_{2}\end{array}\).$$
So then 
\begin{align*}
&\(\begin{array}{c}\(\partial_{u_{1}}F_{1}\(0,0,0\)+1\)\(\Delta_g w_{1}+w_{1}\)\\\(\partial_{u_{2}}F_{2}\(0,0,0\)+1\)\(\Delta_g w_{2}+w_{2}\)\end{array}\)\\
&\qquad=-\(\begin{array}{c}\partial_{u_{1}}\partial_{\alpha}F_{1}\(0,0,0\)\(\Delta_g v_{1}+v_{1}\)\\\partial_{u_{2}}\partial_{\alpha}F_{2}\(0,0,0\)\(\Delta_g v_{2}+v_{2}\)\end{array}\),
\end{align*}
which gives 
$$\(\begin{array}{c}\(\partial_{u_{1}}F_{1}\(0,0,0\)+1\)w_{1}\\\(\partial_{u_{2}}F_{2}\(0,0,0\)+1\)w_{2}\end{array}\)=-\(\begin{array}{c}\partial_{u_{1}}\partial_{\alpha}F_{1}\(0,0,0\)v_{1}\\\partial_{u_{2}}\partial_{\alpha}F_{2}\(0,0,0\)v_{2}\end{array}\).$$
Now, if we suppose that $\(w_{1},w_{2}\)\in\range\(D_{U}\mathcal{T}\(0,0,0\)\)$, then by the assumption (C5), we obtain that $\(v_{1},v_{2}\)\in\range\(D_{U}\mathcal{T}\(0,0,0\)\)$ and so there exists $\(\varphi_1,\varphi_2\)\in\mathcal{H}$ such that 
$$\(\begin{array}{c}\Delta_{g}\varphi_{1}-\partial_{u_{1}}F_{1}\(0,0,0\)\varphi_{1}\\\Delta_{g}\varphi_{2}-\partial_{u_{2}}F_{2}\(0,0,0\)\varphi_{2}\end{array}\)=\(\begin{array}{c}\Delta_g v_{1}+v_{1}\\\Delta_g v_{2}+v_{2}\end{array}\).$$
For $i=1,2$, straightforward integrations by parts then yield
\begin{align*}
\int_M\(\left|\nabla v_i\right|^2+v_i^2\)dv_g&=\int_Mv_i\(\Delta_g\varphi_i-\partial_{u_{i}}F_{i}\(0,0,0\)\varphi_{i}\)dv_g\\
&=\int_M\(\Delta_gv_i-\partial_{u_{i}}F_{i}\(0,0,0\)v_{i}\)\varphi_idv_g=0,
\end{align*}
where $dv_g$ is the volume element with respect to the metric $g$. It follows that $\(v_{1},v_{2}\)=\(0,0\)$. Hence condition \eqref{odd_bif} is satisfied.

\smallskip\noindent {\it Step 3.}  
Finally we prove the expansion \eqref{bif:expan}. 

\smallskip
Let $\(\(\alpha_{m},u_{1,m},u_{2,m}\)\)_{m\in\N}$ be a sequence of solutions to \eqref{DefEq}, such that $\(\alpha_m,u_{1,m},u_{2,m}\)\in I\times\(\mathcal{U}\backslash\left\{\(0,0\)\right\}\)$ and $\(\alpha_m,u_{1,m},u_{2,m}\)\to\(0,0,0\)$ in $I\times C^2\(M\)^2$ as $m\to\infty$. For $i=1,2$, consider the sequence 
$$w_{i,m}:=\varepsilon_m^{-1}u_{i,m },\quad\text{where}\quad\varepsilon_m:=\max\(\left\|u_{1,m}\right\|_{C^{1,\theta}},\left\|u_{2,m}\right\|_{C^{1,\theta}}\)$$
so that 
$$\max\(\left\|w_{1,m}\right\|_{C^{1,\theta}},\left\|w_{2,m}\right\|_{C^{1,\theta}}\)=1.$$ 
Since $\(\alpha_{m},u_{1,m},u_{2,m}\)$ satisfies \eqref{DefEq}, it follows by our assumptions on $F$ that
\begin{align*}
\Delta_{g}w_{i,m}&=\varepsilon_m^{-1}F_{i}\(\alpha_{m},u_{1,m},u_{2,m}\)=\partial_{u_{i}}F_{i}(\alpha_{m},0,0)w_{i,m}+\smallo\(1\)
\end{align*}
uniformly in $M$. Then, by standard elliptic theory, it follows that up to a subsequence $w_{i,m}\to\varphi_{i}$ in $\mathcal{H}_i$ for some function $\varphi_i\in\mathcal{H}_{i}$ satisfying
$$\Delta_{g} \varphi_{i}=\partial_{u_{i}}F_{i}\(0,0,0\)\varphi_{i}\quad\hbox{in } M.$$
Hence $\varphi_{i}$ belongs to $\mathcal{H}_{i}^{*}$ and further $\max\(\left\|\varphi_1\right\|_{C^{1,\theta}},\left\|\varphi_2\right\|_{C^{1,\theta}}\)=1$. It follows that $\(\varphi_1,\varphi_2\)\in\mathcal{H}^*\backslash\left\{\(0,0\)\right\}$ and
$$u_{i,m}=\varepsilon_m w_{i,m}=\varepsilon_m\(\varphi_{i}+\smallo\(1\)\)\quad\hbox{in }\mathcal{H}_i.$$
This completes the proof of Theorem~\ref{Th3}.
\endproof

We can now prove Theorem~\ref{Th1}~(iii) and Theorem~\ref{Th2} by using Theorem~\ref{Th3}. We start with proving Theorem~\ref{Th2}.

\proof[Proof of Theorem~\ref{Th2}]
First note that the system in \eqref{Th2Eq} can be rewritten as
$$\(\begin{array}{c}\Delta_{g}u_{1}\\\Delta_{g} u_{2}\end{array}\)=\(\begin{array}{c}{\tilde{F}}_{1}\(\alpha,u_{1},u_{2}\)\\{\tilde{F}}_{2}\(\alpha,u_{1},u_{2}\)\end{array}\),$$
where 
\begin{align*}
{\tilde{F}}\(\alpha,u_{1},u_{2}\)&=\(\begin{array}{c}{\tilde{F}}_{1}\(\alpha,u_{1},u_{2}\)\\{\tilde{F}}_{2}\(\alpha,u_{1},u_{2}\)\end{array}\)\\
&:=\(\begin{array}{c}a_{11}\(\alpha\)\left|u_{1}\right|^{q-2}u_{1}+a_{12}\(\alpha\)\left|u_{2}\right|^{q-2}u_{1}-\lambda_{1}\(\alpha\)u_{1}\\a_{21}\(\alpha\)\left|u_{1}\right|^{q-2}u_{2}+a_{22}\(\alpha\)\left|u_{2}\right|^{q-2} u_{2}-\lambda_{2}\(\alpha\)u_{2}\end{array}\).
\end{align*}
Now let's transform this system so as to apply Theorem~\ref{Th3}. For every $\alpha\in I$, the unique constant solution to \eqref{Th2Eq} (which existence follows from the assumption (B1)) is given by
$$\left\{\begin{aligned}&\overline{u}_{1}\(\alpha\)=\(\frac{\lambda_{1}\(\alpha\)a_{22}\(\alpha\)-\lambda_{2}\(\alpha\)a_{12}\(\alpha\)}{a_{11}\(\alpha\)a_{22}\(\alpha\)-a_{21}\(\alpha\)a_{12}\(\alpha\)} \)^{1/\(q-2\)}\\
&\overline{u}_{2}\(\alpha\)=\(\frac{\lambda_{2}\(\alpha\)a_{11}\(\alpha\)-\lambda_{1}\(\alpha\)a_{21}\(\alpha\)}{a_{11}\(\alpha\)a_{22}\(\alpha\)-a_{21}\(\alpha\)a_{12}\(\alpha\)}\)^{1/\(q-2\)}\end{aligned}\right.$$
and it satisfies
$$\left\{\begin{aligned}
&a_{11}\(\alpha\)\overline{u}_{1}\(\alpha\)^{q-2}+a_{12}\(\alpha\)\overline{u}_{2}\(\alpha\)^{q-2}=\lambda_{1}\(\alpha\)&&\hbox{ in } M\\
&a_{21}\(\alpha\)\overline{u}_{1}\(\alpha\)^{q-2}+a_{22}\(\alpha\)\overline{u}_{2}\(\alpha\)^{q-2}=\lambda_{2}\(\alpha\)&&\hbox{ in } M. 
\end{aligned}\right.$$
We look for solutions of \eqref{IntroEq1} bifurcating from $\(\overline{u}_{1}\(\alpha\),\overline{u}_{2}\(\alpha\)\)$. By the assumption (B2), for every $\alpha\in I$, the matrix $\mathcal{A}\(\alpha\)$ has two distinct, non-zero, real eigenvalues $\beta_{1}\(\alpha\)$ and $\beta_{2}\(\alpha\)$ given by
$$\left\{\beta_{1}\(\alpha\),\beta_{2}\(\alpha\)\right\}:=\left\{\frac{a_{11}\(\alpha\)\overline{u}_{1}\(\alpha\)^{q-2}+a_{22}\(\alpha\)\overline{u}_{2}\(\alpha\)^{q-2}}{2}\pm\frac{\sqrt{D\(\alpha\)}}{2}\right\},$$
where
\begin{multline*}
D\(\alpha\):=\(a_{11}\(\alpha\)\overline{u}_{1}\(\alpha\)^{q-2}-a_{22}\(\alpha\)\overline{u}_{2}\(\alpha\)^{q-2}\)^{2}\\
+4a_{12}\(\alpha\)a_{21}\(\alpha\)\overline{u}_{1}\(\alpha\)^{q-2}\overline{u}_{2}\(\alpha\)^{q-2}.
\end{multline*}
Let $\mathcal{P}\(\alpha\)$ be the $2\times2$ matrix such that 
$$\mathcal{A}\(\alpha\)=\mathcal{P}\(\alpha\)^{-1}\(\begin{array}{cc}\beta_{1}\(\alpha\)&0\\0&\beta_{2}\(\alpha\)\end{array}\)\mathcal{P}\(\alpha\).$$
Consider $\(u_{1},u_{2}\)\in C^{2}\(M\)^2$ and  let 
$$\(\begin{array}{c}v_{1}\\v_{2}\end{array}\):=P\(\alpha\)\(\begin{array}{c}u_{1}-\overline{u}_{1}\(\alpha\)\\ u_{2}-\overline{u}_{2}\(\alpha\)\end{array}\).$$
We then define
\begin{equation}\label{def:F}
F\(\alpha,v_{1},v_{2}\)=\(\begin{array}{c}F_{1}\(\alpha,v_{1},v_{2}\)\\F_{2}\(\alpha,v_{1},v_{2}\) \end{array}\):=P\(\alpha\)\tilde{F}\(\alpha,u_{1},u_{2}\),
\end{equation}
where
\begin{equation}\label{def:u}
\(\begin{array}{c}u_{1}\\u_{2}\end{array}\)=\(\begin{array}{c}\overline{u}_{1}\(\alpha\)\\ \overline{u}_{2}\(\alpha\)\end{array}\)+P\(\alpha\)^{-1}\(\begin{array}{c}v_{1}\\v_{2}\end{array}\).
\end{equation}
We then obtain that the system \eqref{IntroEq1} is equivalent to
\begin{align}\label{sys:reduced}
\(\begin{array}{c}\Delta_{g}v_{1}\\\Delta_{g}v_{2}\end{array}\)=\(\begin{array}{c}F_{1}\(\alpha,v_{1},v_{2}\)\\F_{2}\(\alpha,v_{1},v_{2}\)\end{array}\).
\end{align}

\smallskip
Next we apply Theorem~\ref{Th3} to \eqref{sys:reduced}. Note that the condition (C1) of Theorem~\ref{Th3} is satisfied by \eqref{sys:reduced}. Furthermore, since $\(\overline{u}_{1}\(\alpha\),\overline{u}_{2}\(\alpha\)\)\in  \(0,\infty\)^2$ for all $\alpha\in I$, by continuity, we obtain that there exists $\delta_0>0$ such that $\(u_1,u_2\)\in \(0,\infty\)^2$ for all $\(\alpha,v_1,v_2\)\in I\times\(-\delta_0,\delta_0\)^2$. In particular, letting $\Omega:=\(-\delta_0,\delta_0\)^2$, we then obtain that $F\in C^1\(I\times\Omega\)$, $\partial_{u_{i}}\partial_\alpha F_{j}$ and $\partial_{u_{i}}\partial_{u_{j}}F_{k}$ exist and are continuous in $I\times\Omega$ for all $i,j,k\in\left\{1,2\right\}$ and the condition (C3) is satisfied with $\mathcal{H}_{1}=\mathcal{H}_{2}:=C^{1,\theta}\(M\)$ and $\mathcal{U}_{1}=\mathcal{U}_{2}:=C^{1,\theta}\(M,\(-\delta_0,\delta_0\)\)$. The condition (C2) is also satisfied as we obtain differentiating  
\begin{align*}
D_{\(v_{1},v_{2}\)}F\(\alpha,0,0\)&=\mathcal{P}\(\alpha\)\[\(q-2\)\mathcal{A}\(\alpha\)\]\mathcal{P}^{-1}\(\alpha\)\\
&=\(q-2\)\(\begin{array}{cc}\beta_{1}\(\alpha\)&0\\0&\beta_{2}\(\alpha\)\end{array}\).
\end{align*}
The assumptions (B3) and (B4) then imply that the conditions (C4) and (C5) of Theorem~\ref{Th3} are also satisfied. 

\smallskip
By applying Theorem~\ref{Th3} and reversing the above change of function, we then obtain that the solution $\(0,\overline{u}_1\(0\),\overline{u}_2\(0\)\)$ is a bifurcation point of the system \eqref{Th2Eq}. Furthermore, we obtain that for every sequence $\(\(\alpha_{m},u_{1,m},u_{2,m}\)\)_{m\in\N}$ of solutions to \eqref{Th2Eq}, if $\(u_{1,m},u_{2,m}\)\ne\(\overline{u}_1\(0\),\overline{u}_2\(0\)\)$ and $\(\alpha_m,u_{1,m},u_{2,m}\)\to\(0,\overline{u}_1\(0\),\overline{u}_2\(0\)\)$ in $I\times C^2\(M\)^2$ as $m\to\infty$, then up to a subsequence, 
$$u_{i,m}=\overline{u}_{i}\(\alpha_{m}\)+\varepsilon_{m}\(q_{i1}\varphi_{1}+q_{i2}\varphi_{2}+\smallo\(1\)\)\quad\text{ in }\mathcal{H}_i\quad\forall i\in\left\{1,2\right\},$$
where $\(q_{ij}\)_{1\le i,j\le 2}:=P\(0\)^{-1}$, for some $\(\varphi_{1},\varphi_{2}\)\in\mathcal{H}^*\backslash\left\{\(0,0\)\right\}$ and $\varepsilon_{m}>0$ such that $\varepsilon_m\to0$. By the assumptions (B2) and (B4), we have that for $i=1,2$, either $\varphi_i\equiv0$ or $\varphi_i$ is not constant in $M$. Also by  assumption (B4) we have that if $\varphi_i\ne0$  for $i=1,2$, then either $\lambda_1\(0\)\ne\lambda_2\(0\)$ or $\beta_{i}\(0\)\ne\lambda_1\(0\)=\lambda_2\(0\)$, which implies $\overline{u}_{2}\(0\)q_{1i}\neq \overline{u}_{1}\(0\)q_{2i}$. In particular, we obtain that $\(u_{1,m},u_{2,m}\)$ is non-synchronized. Therefore, we obtain that there exists a neighborhood $\mathcal{N}$ of $\(0,\overline{u}_{1}\(0\), \overline{u}_{2}\(0\)\)$ in $I\times C^2\(M\)^2$ such that for every solution $\(\alpha,u_1, u_2\)\in\mathcal{N}$ of \eqref{Th2Eq}, if $\(u_{1},u_{2}\)\ne\(\overline{u}_{1}\(\alpha\), \overline{u}_{2}\(\alpha\)\)$, then $\(u_{1},u_{2}\)$ is non-synchronized. If moreover $\mathcal{H}^{*}$ has dimension one, then it follows from the last part of Theorem~\ref{Th3} that there exists a $C^1$ branch of non-synchronized solutions to \eqref{Th2Eq} emanating from $\(0,\overline{u}_{1}\(0\), \overline{u}_{2}\(0\)\)$. This completes the proof of Theorem~\ref{Th2}.
\endproof

\proof[Proof of Theorem~\ref{Th1}~(iii)]
We proceed as in the proof of Theorem~\ref{Th2}. By the assumptions (A1) and (A3) along with the continuity of $\lambda$, $a$ and $b$, letting $\delta$ be smaller if necessary, we may assume that $a\(\alpha\)>b\(\alpha\)$ and $\lambda\(\alpha\)\(a\(\alpha\)+b\(\alpha\)\)>0$ for all $\alpha\in I$. Then the unique constant solution for the system  \eqref{IntroEq4} is given by  
$$\(\begin{array}{c}\overline{u}_{1}\(\alpha\)\\\overline{u}_{2}\(\alpha\)\end{array}\):=\overline{u}\(\alpha\)\(\begin{array}{c}1\\1\end{array}\),\quad\text{where }\overline{u}\(\alpha\):=\(\frac{\lambda\(\alpha\)}{a\(\alpha\)+b\(\alpha\)}\)^{{1}/{\(q-2\)}}$$
and it satisfies
$$\left\{\begin{aligned}&a\(\alpha\)\overline{u}_{1}\(\alpha\)^{q-2}+b\(\alpha\)\overline{u}_{2}\(\alpha\)^{q-2}=\lambda\(\alpha\)\\&b\(\alpha\)\overline{u}_{1}\(\alpha\)^{q-2}+a\(\alpha\)\overline{u}_{2}\(\alpha\)^{q-2}=\lambda\(\alpha\).\end{aligned}\right.$$
We look for solutions of \eqref{IntroEq4} bifurcating from $\(\overline{u}_{1}\(\alpha\), \overline{u}_{2}\(\alpha\)\)$. For \eqref{IntroEq4}, the eigenvalues $\beta_{1}\(\alpha\)$ and $\beta_{2}\(\alpha\)$ of $\mathcal{A}\(\alpha\)$ are given by
$$\beta_{1}\(\alpha\)=\lambda\(\alpha\)\quad\text{and}\quad\beta_{2}\(\alpha\)=\lambda\(\alpha\)\frac{a\(\alpha\)-b\(\alpha\)}{a\(\alpha\)+b\(\alpha\)}\,.$$
We let $F$ and $\Omega=\(-\delta_0,\delta_0\)^2$ be defined similarly as in the proof of Theorem~\ref{Th2}, so that in particular, the conditions (C1) and (C2) of Theorem~\ref{Th3} are satisfied. In this case, we find
$$\mathcal{P}\(\alpha\)=\(\begin{array}{cc}1&1\\1&-1\end{array}\).$$

\smallskip
Next, we choose the appropriate subspaces $\mathcal{H}_{1}$ and $\mathcal{H}_{2}$ and open subsets $\mathcal{U}_1\subseteq\mathcal{H}_1$ and $\mathcal{U}_2\subseteq\mathcal{H}_2$. For this, we use an idea from Gladiali, Grossi and Troestler~\cite{GGT1}. Consider the reflexion $\hat{v}$ across the equator $\left\{x_{n}=0\right\}$ of the sphere $\S^n$ defined by
$$\hat{v}\(x\):=v\(x_{1},\ldots,x_{n},-x_{n+1}\)\quad\forall x=\(x_{1},\ldots,x_{n},x_{n+1}\)\in\S^n$$
for all functions $v:\S^n\to\R$. By stereographic projection along with a conformal change of metric, this corresponds to the Kelvin transform in $\R^{n}$. We let $N_0:=\(0,\ldots,0,1\)$ and $j_0\in\N$ be such that 
$$\(q-2\)\beta_2\(0\)=\lambda_{j_0}:=j_0\(j_0+n-1\)$$ 
i.e. the $j_0$-th eigenvalue of $\Delta_{g_0}$ on $\S^n$ (the existence of $j_0$ is given by the assumption (A3)). We then define $\mathcal{H}_{1}$ and $\mathcal{H}_{2}$ as
$$\mathcal{H}_{1}:=\left\{v\in C^{1,\theta}\(\S^n\):\hbox{$v$ is radial with respect to $N_0$ and $\hat{v}=v$}\right\}$$
and
$$\mathcal{H}_{2}:=\left\{\hspace{-4pt}\begin{aligned}&\left\{v\in C^{1,\theta}\(\S^n\):\hbox{$v$ is radial w.r.t.~$N_0$ and $\hat{v}=v$}\right\}\,\hbox{if $j_0$ is even}\\&\left\{v\in C^{1,\theta}\(\S^n\):\hbox{$v$ is radial w.r.t.~$N_0$ and $\hat{v}=-v$}\right\}\,\hbox{if $j_0$ is odd}.
\end{aligned}\right.$$
For $i=1,2$, we take $\mathcal{U}_{i}:=C^{1,\theta}\(\S^n,\(-\delta_0,\delta_0\)\)\cap \mathcal{H}_i$. For every $\(v_1,v_2\)\in\mathcal{U}$, letting $\(u_1,u_2\)$ be as in \eqref{def:u}, we then obtain
\begin{align*}
\(\begin{array}{c}\hat{u}_{1}\\\hat{u}_{2}\end{array}\)&=\(\begin{array}{c}\overline{u}_{1}\(\alpha\)\\ \overline{u}_{2}\(\alpha\)\end{array}\)+\frac{1}{2}\(\begin{array}{c}\hat{v}_{1}+\hat{v}_{2}\\\hat{v}_{1}-\hat{v}_{2}\end{array}\)\\
&=\left\{\begin{aligned}&\(\begin{array}{c}\overline{u}_{1}\(\alpha\)\\\overline{u}_{2}\(\alpha\)\end{array}\)+\frac{1}{2}\(\begin{array}{c}v_{1}+v_{2}\\v_{1}-v_{2}\end{array}\)=\(\begin{array}{c}u_{1}\\u_{2}\end{array}\)&&\text{if }j_0\text{ is even}\\&\(\begin{array}{c}\overline{u}_{1}\(\alpha\)\\\overline{u}_{2}\(\alpha\)\end{array}\)+\frac{1}{2}\(\begin{array}{c}v_{1}-v_{2}\\v_{1}+v_{2}\end{array}\)=\(\begin{array}{c}u_{2}\\u_{1}\end{array}\)&&\text{if }j_0\text{ is odd}.\end{aligned}\right.
\end{align*}
It follows that
\begin{align*}
F\(\alpha,\hat{v}_{1},\hat{v}_{2}\)&=\(\begin{array}{c}F_{1}\(\alpha,\hat{v}_{1},\hat{v}_{2}\)\\F_{2}\(\alpha,\hat{v}_{1},\hat{v}_{2}\) \end{array}\)\\
&=\(\begin{array}{cc}1&1\\1&-1\end{array}\)\(\begin{array}{c}a\(\alpha\)\hat{u}_{1}^{q-1}+b\(\alpha\)\hat{u}_{2}^{q-2}\hat{u}_{1}-\lambda\(\alpha\) \hat{u}_{1}\\b\(\alpha\)\hat{u}_{1}^{q-2}\hat{u}_{2}+a\(\alpha\)\hat{u}_{2}^{q-1}-\lambda\(\alpha\) \hat{u}_{2}\end{array}\)\\
&=\left\{\begin{aligned}&F\(\alpha,v_{1},v_{2}\)&&\text{if }j_0\text{ is even}\\&\(\begin{array}{c}F_1\(\alpha,v_{1},v_{2}\)\\-F_2\(\alpha,v_{1},v_{2}\)\end{array}\)&&\text{if }j_0\text{ is odd}.\end{aligned}\right.
\end{align*}
This along with standard elliptic regularity and symmetry arguments gives that the condition (C3) is satisfied. 

\smallskip
Recall that the spherical harmonics $\varphi$ satisfying $\Delta_{g_0}\varphi=\lambda_{j_0}\varphi$ are given by the restriction to $\S^n$ of the harmonic polynomials of degree $j_0$ in $\R^{n+1}$. In particular, up to a constant factor, the unique such function $\varphi_{N_0,j_0}$ that is radial with respect to $N_0$ is given by the Jacobi polynomial
\begin{multline*}
\varphi_{N_0,j_0}\(x_1,\dotsc,x_{n+1}\):=\sum_{j=0}^{j_0}\(\begin{array}{c}j_0+\(n-2\)/2\\j\end{array}\)\(\begin{array}{c}j_0+\(n-2\)/2\\j_0-j\end{array}\)\\
\times\(\frac{x_{n+1}-1}{2}\)^{j_0-j}\(\frac{x_{n+1}+1}{2}\)^{j}\quad\forall x\in\S^n
\end{multline*}
(see Gladiali~\cite{G}). So then the assumptions (A2) and (A3) give that $\mathcal{H}^{*}_{1}=\left\{0\right\}$ and $\mathcal{H}^{*}_{2}$ has dimension one. In particular, we obtain that the condition (C4) is satisfied. Furthermore, the condition (C5) follows from (A3).

\smallskip
By applying Theorem~\ref{Th3}, we then obtain that there exists a $C^1$ branch of solutions to \eqref{Th2Eq} emanating from $\(0,\overline{u}\(0\),\overline{u}\(0\)\)$ whose tangent line at $\(0,\overline{u}\(0\),\overline{u}\(0\)\)$ is directed by some vector of the form $\(\mu,\varphi,-\varphi\)$, where $\mu\in\R$ and $\varphi\in \mathcal{H}_2^*\backslash\left\{0\right\}$. In particular, we obtain that near $\(0,\overline{u}\(0\),\overline{u}\(0\)\)$, the solutions on this branch are non-synchronized, which completes the proof of Theorem~\ref{Th1}~(iii).
\endproof

\section{Synchronization and non-existence results}\label{Proportional}

In this section, we prove the following results, which extend Theorem~\ref{Th1} to more general systems and manifolds:

\begin{theorem}\label{Th4} (Case $\lambda_1=\lambda_2$)
Let $\lambda,a_{11},a_{12},a_{21},a_{22}\in\R$, $q\in\(2,\infty\)$ and $\(M,g\)$ be a smooth, closed Riemannian manifold. Consider the system \eqref{IntroEq1} with $\lambda_1=\lambda_2=\lambda$.  
\begin{enumerate}
\item[(i)] If either [$a_{21}\le a_{11}$ and $a_{22}\le a_{12}$] or [$a_{11}\le a_{21}$ and $a_{12}\le a_{22}$] and at least one of the two inequalities is strict, then \eqref{IntroEq1} has no solutions.  
\item[(ii)] If either [$a_{11}<a_{21}$ and $a_{22}<a_{12}$] or [$a_{11}=a_{21}$ and $a_{22}=a_{12}$], then every solution of \eqref{IntroEq1} is synchronized.
\item[(iii)] Assuming that there exists at least one non-zero eigenvalue of $\Delta_g$ with odd multiplicity, Theorem~\ref{Th2} provides examples of real numbers $\lambda,a_{11},a_{12},a_{21},a_{22}>0$ such that $a_{21}<a_{11}$, $a_{12}<a_{22}$ and \eqref{IntroEq1} has non-synchronized solutions.
\end{enumerate}
\end{theorem}

\begin{theorem}\label{Th5} (Case $\lambda_1<\lambda_2$)
Let $\lambda_1,\lambda_2,a_{11},a_{12},a_{21},a_{22}\in\R$, $q\in\(2,\infty\)$ and  $\(M,g\)$ be a smooth, closed Riemannian manifold. Consider the system \eqref{IntroEq1}.   Assume that $\lambda_1<\lambda_2$.
\begin{enumerate}
\item[(i)] If [$a_{21}\le a_{11}$ and $a_{22}\le a_{12}$], then \eqref{IntroEq1} has no solutions.  
\item[(ii)] Assuming that there exists at least one non-zero eigenvalue of $\Delta_g$ with odd multiplicity, Theorem~\ref{Th2} provides examples of real numbers $\lambda_1,\lambda_2,a_{11},a_{12},a_{21},a_{22}>0$ such that $\lambda_2>\lambda_1>0$ and \eqref{IntroEq1} has non-synchronized solutions in each of the following cases: [$a_{21}<a_{11}$ and $a_{12}<a_{22}$],  [$a_{11}<a_{21}$ and $a_{22}<a_{12}$], [$a_{11}<a_{21}$ and $a_{12}<a_{22}$], [$a_{11}=a_{21}$ and $a_{12}<a_{22}$] and [$a_{11}<a_{21}$ and $a_{12}=a_{22}$].
\end{enumerate}
\end{theorem}

\proof[Proof of Theorem~\ref{Th4}~(i) and Theorem~\ref{Th5}~(i)]
Suppose that there exists a solution $\(u_1,u_2\)$ of \eqref{IntroEq1}. We define 
$$v\(x\):=\frac{u_1\(x\)}{u_2\(x\)}\qquad\forall x\in M.$$ 
We then obtain 
\begin{align*}
&\Delta_{g}v=\frac{\Delta_{g}u_{1}}{u_{2}}+2\frac{\<\nabla_{g} u_{1},\nabla_{g}u_{2}\>_{g}}{u_{2}^{2}}-\frac{u_{1}\Delta_{g}u_{2}}{u_{2}^{2}}-2\frac{u_{1}\left|\nabla_{g}u_{2}\right|^{2}_{g}}{u_{2}^{3}}\\
&=\(a_{11}-a_{21}\)u_{1}^{q-2}v+\(a_{12}-a_{22}\)u_{2}^{q-2}v+\(\lambda_{2}-\lambda_{1}\)v+2\frac{\<\nabla_{g}v,\nabla_{g}u_{2}\>_{g}}{u_{2}}
\end{align*}
in $M$. In the case where $a_{21}\le a_{11}$, $a_{22}\le a_{12}$, $\lambda_1\le\lambda_2$ and one of these inequalities is strict, we then obtain 
$$\Delta_{g}v>2\<\nabla_{g}v,\nabla_{g}\[\ln u_2\]\>_{g}\quad\text{in }M,$$
which is in contradiction with the minimum principle. Similarly, in the case where $a_{11}\le a_{21}$, $a_{12}\le a_{22}$, $\lambda_2\le\lambda_1$ and one of these inequalities is strict, we obtain 
$$\Delta_{g}v<2\<\nabla_{g}v,\nabla_{g}\[\ln u_2\]\>_{g}\quad\text{in }M,$$
which is in contradiction with the maximum principle. This proves Theorem~\ref{Th4}~(i) and Theorem~\ref{Th5}~(i).
\endproof

\proof[Proof of Theorem~\ref{Th4}~(ii)]
Suppose first that $a_{11}<a_{21}$, $a_{22}<a_{12}$ and $\lambda_{1}=\lambda _{2}$. Let $v$ be as in the previous proof and $x_{1},x_{2}\in M$ be such that 
$$v\(x_{1}\)=\min\left\{v\(x\):\ x\in M\right\}\quad\text{and}\quad v\(x_{2}\)=\max\left\{v\(x\):\ x\in M\right\}.$$ 
We then obtain 
$$0\ge\Delta_{g}v\(x_{1}\)=\[\(a_{12}-a_{22}\)u_{2}\(x_{1}\)^{q-2}-\(a_{21}-a_{11}\)u_{1}\(x_{1}\)^{q-2}\]v\(x_{1}\)$$
and
$$0\le\Delta_{g}v\(x_{2}\)=\[\(a_{12}-a_{22}\)u_{2}\(x_{2}\)^{q-2}-\(a_{21}-a_{11}\)u_{1}\(x_{2}\)^{q-2}\]v\(x_{2}\),$$
which imply
$$v\(x_{2}\)\le\(\frac{a_{12}-a_{22}}{a_{21}-a_{11}} \)^{1/\(q-2\)}\le v\(x_{1}\).$$
It follows that $v$ is constant in $M$. 

\smallskip
Now suppose that $a_{11}=a_{21}$, $a_{22}=a_{12}$ and $\lambda_{1}=\lambda _{2}$. In this case, we have 
$$\Delta_{g}v=2\<\nabla_{g}v,\nabla_{g}\[\ln u_2\]\>_{g}\quad\text{in }M.$$
It then follows from the maximum principle that $v$ is constant in $M$. This completes the proof of Theorem~\ref{Th4}~(ii).
\endproof

\proof[Proof of Theorem~\ref{Th4}~(iii) and Theorem~\ref{Th5}~(ii)]
We choose our examples of the form $a_{11}\(\alpha\)=\lambda_{1}\(\alpha\)a$, $a_{12}\(\alpha\)=\lambda_{1}\(\alpha\)b$, $a_{21}\(\alpha\)=\lambda_{2}\(\alpha\)b$ and $a_{22}\(\alpha\)=\lambda_{2}\(\alpha\)a$ for some $a,b>0$ and $\lambda_1,\lambda_2\in C^1\(I,\(0,\infty\)\)$, where $I:=\[-\delta,\delta\]$ for some $\delta>0$ to be chosen small. As is easily seen, for every $\alpha\in I$, if 
$$a\ne b\quad\text{and}\quad D\(\alpha\):=\(\lambda_{1}\(\alpha\)-\lambda_{2}\(\alpha\)\)^{2} a^{2}+4\lambda_{1}\(\alpha\)\lambda_{2}\(\alpha\)b^{2}>0,$$
then the system \eqref{IntroEq1} has a unique constant solution given by
$$\(\begin{array}{c}\overline{u}_{1}\(\alpha\)\\\overline{u}_{2}\(\alpha\)\end{array}\)=\(a+b\)^{-1/\(q-2\)}\(\begin{array}{c}1\\1\end{array}\).$$
and the matrix $A\(\alpha\)$ has two distinct, non-zero real eigenvalues given by
$$\left\{\beta_{1}\(\alpha\),\beta_{2}\(\alpha\)\right\}=\left\{\frac{\(\lambda_{1}\(\alpha\)+\lambda_{2}\(\alpha\)\)a\pm\sqrt{D\(\alpha\)}}{2\(a+b\)}\right\}.$$
Now, we treat each case separately and choose $a,b,\lambda_1$ and $\lambda_2$ in such a way that 
$a\ne b$, $D>0$, $\(q-2\)\beta_1\(0\)\not\in\spec\(\Delta_g\)$, $\(q-2\)\beta_2\(0\)=\lambda_0$, $\beta'_2\(0\)\ne0$ and either $\lambda_1\(0\)\ne\lambda_2\(0\)$ or $\beta_2\(0\)\ne\lambda_1\(0\)=\lambda_2\(0\)$, where $\spec\(\Delta_g\)$ stands for the spectrum of $\Delta_g$ and $\lambda_0$ is a non-zero eigenvalue of $\Delta_g$ with odd multiplicity (which existence is given by assumption). Choosing $\delta$ small enough, we can then apply Theorem~\ref{Th2}.

\smallskip\noindent 
{\it Case $\lambda_{1}=\lambda_{2}$, $a_{21}<a_{11}$ and $a_{12}<a_{22}$.}
Choose for instance $a:=\frac{\lambda_{0}+\lambda}{q-2}$, $b:=\frac{\lambda-\lambda_{0}}{q-2}$ and $\lambda_{1}\(\alpha\)=\lambda_{2}\(\alpha\):=\frac{\lambda\(\alpha+1\)}{q-2}$ for some $\lambda\in\(\lambda_0,\infty\)\backslash\spec\(\Delta_{g}\)$. Then $\lambda_1=\lambda_2>0$, $0<a_{21}=a_{12}<\frac{\lambda\lambda_1}{q-2}<a_{11}=a_{22}$, $a\ne b$, $D=4\lambda_1^2b^2>0$, $\beta_{1}\(0\)=\frac{\lambda}{q-2}$ (so that $\(q-2\)\beta_1\(0\)\not\in\spec\(\Delta_g\)$) and $\beta_{2}\(0\)=\beta_{2}'\(0\)=\frac{\lambda_{0}}{q-2}>\frac{\lambda}{q-2}=\lambda_1\(0\)=\lambda_2\(0\)$.

\smallskip\noindent 
{\it Case $\lambda_{1}<\lambda_{2}$, $a_{21}<a_{11}$ and $a_{12}<a_{22}$.} 
Choose for instance $a:=1$, $b:=\varepsilon$, $\lambda_1\(\alpha\):=\frac{2\(1+\varepsilon\)\lambda_0\(\alpha+1\)}{(2+\varepsilon+\varepsilon\sqrt{5+4\varepsilon})\(q-2\)}$ and $\lambda_2\(\alpha\):=\(1+\varepsilon\)\lambda_1\(\alpha\)$ for some small $\varepsilon>0$. Then $0<\lambda_1<\lambda_2$, $0<a_{21}=\varepsilon\(1+\varepsilon\)\lambda_1<\lambda_1=a_{11}$, $0<a_{12}=\varepsilon\lambda_1<\(1+\varepsilon\)\lambda_1=a_{22}$, $a\ne b$, $D=\varepsilon^2\(5+4\varepsilon\)\lambda_1^2>0$, $\beta_{1}\(0\)=\frac{(2+\varepsilon-\varepsilon\sqrt{5+4\varepsilon})\lambda_0}{(2+\varepsilon+\varepsilon\sqrt{5+4\varepsilon})\(q-2\)}<\frac{\lambda_0}{q-2}$, $\beta_{1}\(0\)\to\frac{\lambda_0}{q-2}$ as $\varepsilon\to0$ (so that $\(q-2\)\beta_1\(0\)\not\in\spec\(\Delta_g\)$) and $\beta_{2}\(0\)=\beta_{2}'\(0\)=\frac{\lambda_{0}}{q-2}$.

\smallskip\noindent 
{\it Case $\lambda_{1}<\lambda_{2}$, $a_{11}<a_{21}$ and $a_{22}<a_{12}$.} 
Choose for instance $a:=1$, $b:=\sqrt{6}$, $\lambda_{1}\(\alpha\):=\frac{(1+\sqrt{6})\lambda_{0}\(\alpha+1\)}{5\(q-2\)}$ and $\lambda_{2}\(\alpha\):=2\lambda_{1}\(\alpha\)$. Then $0<\lambda_1<\lambda_2$, $0<a_{11}=\lambda_1<2\sqrt6\lambda_1=a_{21}$, $0<a_{22}=2\lambda_1<\sqrt6\lambda_1=a_{12}$, $a\ne b$, $D=49\lambda_1^2>0$, $\beta_{1}\(0\)=\frac{-2\lambda_{0}}{5\(q-2\)}<0$ (so that $\(q-2\)\beta_1\(0\)\not\in\spec\(\Delta_g\)$) and $\beta_{2}\(0\)=\beta_{2}'\(0\)=\frac{\lambda_{0}}{q-2}$.

\smallskip\noindent 
{\it Case $\lambda_{1}<\lambda_{2}$, $a_{11}<a_{21}$ and $a_{12}<a_{22}$.} 
Choose for instance $a:=1$, $b:=2$, $\lambda_{1}\(\alpha\):=\frac{3\lambda_{0}\(\alpha+1\)}{(3+2\sqrt{6})\(q-2\)}$ and $\lambda_{2}\(\alpha\):=5\lambda_{1}\(\alpha\)$. Then $0<\lambda_1<\lambda_2$, $0<a_{11}=\lambda_1<10\lambda_1=a_{21}$, $0<a_{12}=2\lambda_1<5\lambda_1=a_{22}$, $a\ne b$, $D=96\lambda_1^2>0$, $\beta_{1}\(0\)=\frac{(3-2\sqrt{6})\lambda_{0}}{(3+2\sqrt{6})\(q-2\)}<0$ (so that $\(q-2\)\beta_1\(0\)\not\in\spec\(\Delta_g\)$) and $\beta_{2}\(0\)=\beta_{2}'\(0\)=\frac{\lambda_{0}}{q-2}$.

\smallskip\noindent 
{\it Case $\lambda_{1}<\lambda_{2}$, $a_{11}=a_{21}$ and $a_{12}<a_{22}$.} 
Choose for instance $a:=1+\varepsilon$, $b:=1$, $\lambda_1\(\alpha\):=\frac{2\(2+\varepsilon\)\lambda_0\(\alpha+1\)}{((2+\varepsilon)(1+\varepsilon)+\sqrt{\(1+\varepsilon\)\(4+\varepsilon^2+\varepsilon^3\)})\(q-2\)}$ and $\lambda_2\(\alpha\):=\(1+\varepsilon\)\lambda_1\(\alpha\)$ for some small $\varepsilon>0$. Then $0<\lambda_1<\lambda_2$, $a_{11}=a_{21}=\(1+\varepsilon\)\lambda_1>0$, $0<a_{12}=\lambda_1<\(1+\varepsilon\)^2\lambda_1=a_{22}$, $a\ne b$, $D=\(1+\varepsilon\)\(4+\varepsilon^2+\varepsilon^3\)\lambda_1^2>0$, $\beta_{1}\(0\)=\frac{((2+\varepsilon)(1+\varepsilon)-\sqrt{\(1+\varepsilon\)\(4+\varepsilon^2+\varepsilon^3\)})\lambda_0}{((2+\varepsilon)(1+\varepsilon)+\sqrt{\(1+\varepsilon\)\(4+\varepsilon^2+\varepsilon^3\)})\(q-2\)}>0$, $\beta_{1}\(0\)\to0$ as $\varepsilon\to0$ (so that $\(q-2\)\beta_1\(0\)\not\in\spec\(\Delta_g\)$) and $\beta_{2}\(0\)=\beta_{2}'\(0\)=\frac{\lambda_{0}}{q-2}$.

\smallskip\noindent {\it Case $\lambda_{1}<\lambda_{2}$, $a_{11}<a_{21}$ and $a_{12}=a_{22}$.} 
Choose for instance $a:=1$, $b:=3$, $\lambda_{1}\(\alpha\):=\frac{2\lambda_{0}\(\alpha+1\)}{(1+\sqrt{7})\(q-2\)}$ and $\lambda_{2}\(\alpha\):=3\lambda_{1}\(\alpha\)$. Then $0<\lambda_1<\lambda_2$, $0<a_{11}=\lambda_1<9\lambda_1=a_{21}$, $a_{12}=a_{22}=3\lambda_1>0$, $a\ne b$, $D=112\lambda_1^2>0$, $\beta_{1}\(0\)=\frac{(1-\sqrt{7})\lambda_{0}}{(1+\sqrt{7})\(q-2\)}<0$ (so that $\(q-2\)\beta_1\(0\)\not\in\spec\(\Delta_g\)$) and $\beta_{2}\(0\)=\beta_{2}'\(0\)=\frac{\lambda_{0}}{q-2}$.

\smallskip
The results then follow by applying Theorem~\ref{Th2}.
\endproof

\end{document}